\documentclass[10pt]{amsproc}

\usepackage{amssymb}

\usepackage{graphicx}

\usepackage{comment,color}

\newtheorem{theorem}{Theorem}[section]

\theoremstyle{definition}
\newtheorem{definition}[theorem]{Definition}

\newtheorem{defn}[theorem]{Definition}

\newtheorem{corollary}[theorem]{Corollary}
\theoremstyle{remark}

\usepackage{comment,color}

\newtheorem{rem}[theorem]{Remark}
\numberwithin{equation}{section}

\def\e{\epsilon}

\def\t{\tau}

\def\e{\epsilon}

\def\CC{\mathbb C}

\def\ZZ{\mathbb Z}

\def\RR{\mathbb R}

\def\Tr{\rm Tr }

\def\G_0{\Gamma_0}
\def\G{\Gamma}

\def\t{\tau}

\def\e{\epsilon}

\def\CC{\mathbb C}

\def\ZZ{\mathbb Z}

\def\RR{\mathbb R}

\def\Tr{\rm Tr }
\begin{document}

\title[Computing Adapted Bases]{Computing Adapted Bases for Conformal Automorphism Groups of Riemann Surfaces}

%\maketitle

%    author one information
\author{Jane Gilman}
\address{Mathematics Department, Rutgers University, Newark, NJ 07102}
%\curraddr{}
\email{gilman@rutgers.edu}

\thanks{Some of this work was carried out while the author was a supported visitor at ICERM, the Institute for Computational and Experimental Research in Mathematics}
\subjclass{20,32,30,52}
\keywords{Conformal Automorphism, Riemann surface, Schreier-Reidemeister, Homology Basis, Adapted Basis, Matrix Representation}

\date{revised February 10, 2014}

\dedicatory{This paper is dedicated to Emilio Bujalance on his 60th birthday}%    author two information

\begin{abstract} The concept of an adapted homology basis for a prime order conformal automorphism of a compact Riemann surface originated in \cite{thesis, JG1, Moduli, matrix} and is extended to arbitrary finite groups of conformal automorphisms in \cite{prep}. Here we compute some examples of adapted homology bases for  some groups of automorphisms. The method is to begin by apply the Schreier-Reidemeister rewriting process
and the Schreier-Reidemeister Theorem from \cite{MKS} and then to eliminate generators and relations until there is one single large defining relation for the fundamental group in which every generator and its inverse occurs. We  are then able to compute the action of the group on the  homology image of these generators in the first homology group. The matrix of the action is in a simple form.
\end{abstract}

\maketitle
\section{Introduction}
This paper is about adapted homology bases for conformal automorphism groups  of compact  Riemann surfaces.  The concept of an adapted homology basis can be extended from a cyclic group of prime order to an arbitrary conformal automorphism group \cite{prep}.  Our emphasis here is on the method. Adapted bases can be considered from  two different equivalent approaches: one is from the geometric/topological perspective of coverings and curve lifting and the other is from the algebraic perspective of groups and their subgroups using the Schreier-Reidemeister  theory. In particular, although  in the original concept and initial results
  \cite{thesis, JG1}   used curve lifting, it became clear that using Shreier-Reidemeister theory \cite{MKS} gave a much less ad hoc method and more robust, complete results \cite{matrix, Jalg}. Here we give three explicit examples (sections \ref{section:exone}, \ref{section:extwo}, and  \ref{section:exthree}) demonstrating how the Schreier-Reidemeister method is used to find adapted homology bases.

\section{The Set-Up and Notation}

Let $S$, a compact Riemann surface of genus $g \ge 2$
and
$G$, a group of conformal automorphisms of $S$, with order,  $o(G) =n$.
$G$ acts on $S$ with quotient $S_0=S/G$, of genus $g_0$.

$G$ induces an action on the free homotopy gropup of $S$  and on the first homology of $S$. Here by abuse of notation  we do not distinguish between the action of $G$ on $S$, on the free homotopy group, or on  homology.

\subsection{Well-known facts} \label{section:notationH}

If $y_1,...,y_m$ are the fixed points on $S$ of elements of $G$, then
$G_{y_i}$, the stabilizer of
$y_i$ is cyclic for each $i$. Set $n_i = o(G_{y_i})$, the order of $G_{y_i}$.
The covering $\pi: S \rightarrow S_0$ is branched over $t$ points $p_1,...,p_t$. If $\pi(y_i) = p_j$,
then there are ${\frac{n}{n_i}}$ distinct points projecting onto $p_j$. As this is a regular covering, the points over each given $p_j$ have the same branch order.  The Riemann-Hurwitz Theorem implies
$$ 2g = 2n(g_0) - 2n + 2 + n \Sigma_{i=1}^t( 1- {\frac{1}{n_i}}).$$
We assume that $S$ and $S_0$ are uniformized by the Fuchsian groups $\G$ and $\G_0$ and that $U= \{ x + iy| x, y \in \RR  y>0 \} $ is the upper half-plane in $\CC$ so that $S = U/\Gamma$ and $S_0 = U/\Gamma_0$. While $\G$ is a surface group, $\G_0$ is either a surface group or an orbifold group depending upon whether it contain elliptics. We note that $ \Gamma_0/\Gamma$ is isomorphic to $G$  and that $\exists \phi: \Gamma_0 \rightarrow G$ with $\Gamma  = Ker \phi$.

The map $\phi$ is the known as surface kernel map or
{\sl the surface kernel homomorphism}. A vector that gives the same information as the surface kernel map is known as {\sl the generating vector}. It is the vector
$$(\phi(a_1),...,\phi(a_{g_0}), \phi(b_1),...,\phi(b_{g_0}),\phi(x_1),...,\phi(x_t))$$
where
$$\Gamma_0 = \langle a_1,...,a_{g_0}, b_1,...,b_{g_0}, x_1,...,x_t \; | \;
R=(\Pi_{j=1}^{g_0} [a_j,b_j])\cdot x_1 \cdots x_t = 1; \;  x_i^{n_i}=1, i=1,...,t \rangle .$$
$\Gamma_0/ \Gamma$ is isomorphic to $G$ and acts as outer automorphisms on $\Gamma.$ Both $G$ and $\G_0/\G$ induce an action on the first homology group of $S$. By abuse of notation we do not distinguish between the action of $G$ on the surface and the action induced on the first homology group. We write the fundamental group multiplicatively and the first homology group additively. We use $=$ for equality both in $\G$ and in the homology group, but sometimes  emphasize equality in homology by using $\approx_h$.

\subsection{Original motivation and potential applications}
%and potential applications}	
 The  original motivation for considering adapted bases was to establish an algebraic structure on the moduli space of compact surfaces with punctures.  Moduli space is the quotient of the Teichm\"uller Space under the action of the Mapping-class group, but it can also be obtained by factoring through the Siegel Upper Half plane, that is by taking the images of the period matrices in the Siegel Upper half Plane factored out by the action of  the Siegel Modular group.  It was known by the work of Bailey using the Satake Compactification of the quotient of the Siegel Upper half Plane under the Siegel Modular group, that the moduli space of compact surfaces  was a quasi-projective variety.

 If one begins with a punctured surface, one can form a two sheeted covering that is branched over the punctures. The covering will have a conformal involution that fixes the points lying over the punctures. The matrix of the involution with respect to the adapted basis will then have an image in the Siegel modular group and as such will act on the period matrices of all surfaces. The fixed points then of the involution will be a subvariety corresponding to punctured surfaces and its image in moduli space an algebraic variety.  There is some work to extend the action of the image of the adapted automorphism in the Siegel Modular group to the compactification of moduli space.

We mention a number of potential applications of adapted bases which will be addressed more fully in \cite{prep}.
One area of applications is to the representation variety \cite{wmg1, wmg2, Wentworth}.
We consider the representation variety of the fundamental group of a compact surface of genus $g \ge 2$,
$Hom(\pi_1(S),PSL(2,\RR)/PSL(2,\RR))$ where the action of $PSL(2,\RR)$ is given by conjugation.
  If  $S$ is a surface with a conformal automorphism group $G$, there is a representation of the fundamental group of $S$ where the images of the curves in $PSL(2,\RR)$ have equal traces. That is, there is a set of generators for the fundamental group of $S$ whose images give an adapted homology basis and, therefore, whose traces satisfy certain relations. Conversely, if $\pi_1(S)$ has a representation whose traces satisfy these relations, $S$ has a conformal automorphism group isomorphic to $G$. That is, the length function will have certain symmetries. A corollary of the definition of an adapted basis, definition \ref{defn:prep} below, and the existence of such a basis is:
\begin{corollary} \label{corollary:reps}
Representations of discrete groups with conformal automorphisms can be identified by the traces of elements in the representation.
\end{corollary}

The existence of adapted bases will give further results about strata in the Teichmuller space  and their algebraic structure. It would be nice if it these together with the results for all finite groups could be used to say something new about the disparity in dimensions count for  $\mathcal{H}_g$, the Siegel Upper-half plane (i.e. ($g(g+1)/2$))  and that of the period matrices for compact Riemann surfaces that it contains since Teichmuller space or moduli space has dimension ($3g-3$) for compact surfaces.

\section{Adapted homology basis} \label{section:adapt}
Heuristically, an adapted basis is basis in which the action of a group is in a particularly nice simple form.
The original definition for prime order automorphisms \cite{thesis, JG1,matrix} can be extended to any finite group $G$ of conformal automorphisms.

\begin{definition} %\label{Adapted}
 (\textit{Adapted Homology Basis})  \label{defn:prep}  \cite{prep} The homology basis of $S$,  $\mathcal{B}$,  is {\sl adapted } to a finite group $G$ if for each $\gamma \in \mathcal{B}$ one of the following occurs:
 \begin{enumerate}
 \item \label{item:prop1}  $\gamma$ and $g(\gamma)$ are in the basis for all $g \in G$ and $g(\gamma) \ne  \gamma$.
 \item \label{item:prop2} $\gamma$ and $h^j(\gamma)$ are in the basis for all $j=0, 1,...,m_{i-2}$ where $ h \in G$ is of order $m_i$,  and
         $$h^{m_{i-1}}({\gamma})  =  - (\gamma + h(\gamma) + \cdots  + h^{m_{i-2}}(\gamma)).$$
          Further
       for each right coset representative, $g_h$,  for $G$ modulo $\langle h \rangle$, we have that $g_h(\gamma)$ and $(g_h \circ h^j)(\gamma)$ are in the basis for all $j=0, 1,...,m_{i-2}$ and
       $$(g_h \circ h^{m_{i-1}})({\gamma}) = - (g_h(\gamma) + (g_h \circ h)(\gamma) + \cdots (g_h \circ h^{m_{i-2}})(\gamma)).$$
     \item \label{item:prop3} $\gamma = h^r(\gamma_0)$, $r \ne m_{i-1}$  where $\gamma_0$ is one of the curves in item \ref{item:prop2} above.
        % \item \label{item:prop4} $g(\gamma) = \gamma$ for all $g \in G_0$, $G_0$ a subgroup of $G$ of order $m$.  All of the other $n/m$ images of $\gamma$ under $G$ are fixed appropriately by conjugate elements or  by elements representing the cosets of $G/G_0$ and are also in $\mathcal{B}$.
             \item \label{item:prop4}  $g(\gamma) \in \mathcal{B}$ for all $g \in G_0$, $G_0$ a subgroup of $G$ of order $m$.  All of the other $n/m$ images of $\gamma$ under $G$ are permuted appropriately by conjugate elements. Further if $h$ is a coset representative of $G/G_0$, then either $h(\gamma)$ or $-h(\gamma)$ is in $\mathcal{B}$. \end{enumerate}
 \end{definition}

Here by abuse of notation we do not distinguish between the action of $G$ on the surface and the action induced on the first homology group or between a curve on the surface and its image in homology. We do not use different notation for equality in the group $G$, the fundamental group, or the homology group except we use $\approx_h$ when it seems important to emphasize that this is equality of curves in homology.

When we construct adapted homology bases we refer to curves of type (\ref{item:prop1}), (\ref{item:prop2}), (\ref{item:prop3}), and (\ref{item:prop4}) as in definition \ref{defn:prep} above.

We note that they can be more than one homology basis that is an adapted basis for any given group.

The question of the existence of an adapted homology basis for any $G$ of finite order is addressed in \cite{prep}. Further, the numbers of the different types of elements in the basis depends upon the fixed point structure of the group $G$ which is determined by the surface kernel homomorphism or equivalently by the generating vector in a manner that can be made precise again using definition \ref{defn:prep} and the existence to give

\begin{corollary} \label{corollary:blocks}
\cite{prep} If $G$ is a group of conformal automorphisms of $S$ of order $n$, then the matrix of the action of an element of $G$ on an appropriately ordered adapted homology basis breaks up into blocks: an element of type (\ref{item:prop1}) will give a block that is an  $n \times n$ permutation matrix; elements of type (\ref{item:prop2})  will gives blocks that are $(m_i-1)\times (m_i-1)$  matrices  with $1$'s along the super diagonal, all elements of the last row $-1$ and all other entries $0$; elements of type (\ref{item:prop4}) will give $n \times n$  square matrices that contain permuted blocks of $n/{m_i}$  blocks of $(m_i) \times (m_i)$ permutation  matrices.
\end{corollary}
Note that every element of $G$ either fixes this basis (but not the ordered basis) or maps some elements into their inverse. That is, elements of $G$ will permute the blocks and/or the entries within the blocks, but may map some generators into their inverses.

\begin{rem}
The matrix will, of course, not necessarily be a symplectic matrix, but the intersection numbers of the curves
can be computed from the surface kernel map $\phi$ using the methods of \cite{GP}. The method there is applied to  the prime order case and is due to Patterson, but it also applies to an arbirtary finite order group and has been used by subsequent authors. Further,
 and there are methods for finding the corresponding symplectic matrix and the action of the group on the corresponding homology basis \cite{GP, Jalg}.
\end{rem}
\begin{rem}
We note that if the group has no fixed points, then all elements are of type \ref{item:prop1} or \ref{item:prop4}. If the quotient has genus $0$, all elements are of type \ref{item:prop2} or \ref{item:prop3}, so a large automorphism group in the sense of Kulkarni with quotient of genus $0$, will have all elements $\gamma$ of type \ref{item:prop2} or \ref{item:prop3} and every curve on the surface will satisfy $\Sigma_{g \in G} g(\gamma) \approx_h  1$. Note that Kulkarni defines a large automorphism groups one where the order of the group is greater than $4(g-1)$ so that the genus of the quotient must be $0$ or $1$. The latter occurs only for triangle groups. Here we speak of {\sl large automorphism groups in the sense of Kulkarni of genus $0$}.
\end{rem}
\section{Outline of the Method}
In order to find the adapted basis pick a set of right coset representatives for $\Gamma_0 / \Gamma$. One must choose a Schreier system of cosets, that is one where the initial segment of any representative is itself a representative. Use the theorem from \cite{MKS} Theorem \ref{theorem:MKS} page and given below. This gives a presentation  for the kernel, that is generators and relations for the kernel. This uses the coset representatives to obtain generators for the kernel. There is a rewriting system $\tau$ rewriting any word in the kernel as a word in the chosen generators.
 One  then eliminates of generators and relations via an algorithm \cite{prep} using Teitze transformations until there is only one main relation of length $4g$. The relation contains each word and its inverse exactly once. One can compute the action of elements of the group by conjugating the word by a coset representative on a word using and computing the image of the conjugate under $\tau$.

The full algorithm for the elimination of generators and relations is given in \cite{prep}.
 We give examples of the computations and the full elimination process in the case of three groups.

Treatments of the Schreier-Reidemeister methods can be found in Lyndon-Schupp \cite{LS} and the text of D. Johnson \cite{DaveyJ}. See also works of Marston Condor, Gareth Jones \cite{Jones} and David Singerman as well as that of Costa and Turbek \cite{Costa}.

\section{Schreier-Reidemeister Theorem}

We remind the reader of the Schreier-Reidemeister Theorem.
\begin{theorem}  {\rm(Schreier-Reidemeister \cite{MKS} as stated for  our situation)}
\label{theorem:MKS}
Let $\Gamma_0$ have generators
 $$ a_1,...,a_{g_0}, b_1,...,b_{g_0}, x_1, ...x_t $$
and relations
 $$ R=  \Pi_{i=1}^{g_0} [a_i,b_i]  \cdot x_1 \cdots x_t =1,    \;\;\; x_j^{n_j}=1, j = 1, ..., t.$$
and let $\Gamma$ be the subgroup of $\Gamma_0$ with $G$ as above $G$ isomorphic $\Gamma_0 / \Gamma$.

If $\tau$ is a Reidemeister-Schreier rewriting process, then $\Gamma$ can be presented as a group with generators
$$ S_{K,a_i}, \;\; S_{K,b_i}, i =1,...,g_0 \mbox{ and, } \;\;  S_{K,x_j}, \; j=1,...,t $$
and relations
$$\tau(KRK^{-1}) = 1,\;\; \tau(Kx_j^{n_j}K^{-1})=1,\;\; j = 1,...,t, \;\; S_{M,a}=1$$
where $K$ runs over a complete set of coset representative and $M$ is a coset representative and $a$ any element of $\Gamma_0$ with $Ma$ freely equal to ${\overline{Ma}}.$
\end{theorem}

Here $K$ is a system of Schreier right coset representatives for $\G_0/\G$ (that is, a systems of coset representatives where any initial segment of a representative is also a representative). We let ${\overline{w}}$ denote the coset representative for a word $w \in \G_0$.
$S_{K,y}= Ky{\overline{Ky}}^{-1}$ for any generator $y$. The rewriting process $\tau$ expresses a word in the generators for $\Gamma_0$ as a word in the generators for the kernel, if the word is in the kernel and is defined as follows: let $c_i, i = 1, ..., r$ be generators of $\G_0$ and $\e_i = \pm 1$.
$$\tau(c_1^{\e_1}c_2^{\e_2}c_3^{\e_3} \cdots c_r^{\e_r}) =
S_{K_1,c_1}^{\e_1}S_{K_2,c_2}^{\e_2} \cdots S_{K_r,c_r}^{\e_r}$$
where $K_i= {\overline{c_1c_2\cdots c_{i-1}}}$ if $\e_i=1$ and
$K_i= {\overline{{c_1c_2 \cdots c_{i-1}c_i}}}$ if $\e_i=-1$.

\begin{defn} For ease of reference, when $q$ is in the generating set for $\G_0$, we refer to the $S_{K,q}$ as the $q$-generators for $\G$.
\end{defn}

One can speak of relators for the groups $\G_0$ and $\G$ or of relations. Relators omit the equal sign, but relations do not. We also  use  the equal sign to denote equality  in homology  but  use  $\approx_h$ to denote equality in the homology group for emphasis at times. That is, we use $\approx_h$ to denote homology when the distinction between $=$ and $\approx_h$ is significant.

{\bf Reminder:} As noted in section \ref{section:adapt} we use the same notation for the conjugation action of an element of $\G$ on $\G_0$ as for its induced action on curves in the corresponding homology basis and/or homotopy basis.

\section{Example \#1: A group with no fixed points:  $G = \ZZ_6$} \label{section:exone}

\smallskip

 {\sl The Fuchsian group $\G_0$:} $\G_0= \langle a,b,c,d\;|\; [a,b][c,d]|\rangle$
\vskip .05in
 {\sl The finite group:} $G= \mathbb{Z}_6  = \langle h \rangle$.
\vskip .05in
{\sl The homomorphism, $\phi$:}  $\phi(q)=\phi(b)=
\phi(c) =1; \phi(d)=h.$
{\sl The generating vector:} $(\phi(a),\phi(b),\phi(c),\phi(d)) = (1,1,1,h)$
\vskip .05in
{\sl The Coset representatives:}  $\{1,,d,d^2,d^3,d^4, d^5\}$
\vskip .05in
A  Riemann-Hurwitz calculation gives $g=7$ as $g_0=2$.
\subsection{Subgroup generators}

\smallskip
\subsubsection*{{\bf For $a$}}

Write the subgroup generators for $a$,  $S_{1,a},S_{b,a},S_{d,a},S_{d^2,a},S_{bd,a},S_{bd^2,a}$, as:
\begin{itemize}
\item
$S_{1,a}=a{\overline{a}}^{-1}=: A$
\item

$S_{d,a}=da{\overline{da}}^{-1}= h(A)$
\item
$S_{d^d,a}=d^2a{\overline{d^2}}^{-1} =h^2(A)$
\item
$S_{d^3,a}= d^3a{\overline{d^3}}^{-1} =h^3(A)$
\item
$S_{d^4,a}= d^4a{\overline{d^4a}}^{-1}=h^4(A)$
\item
$S_{d^5,a}= d^5a{\overline{d^5a}}^{-1}= h^5(A)$
\end{itemize}

{\sl Compute} $\t(d^jad^{-j}) = h^j(A)$.
\vskip .05in
\subsubsection*{\bf Computation for images of b and c are similar}
Set $S_{1,b} =B$ and $S_{1,c}=C$ (i.e. $S_{1,b} :=B$ and $S_{1,c}:=C$) so that we have
 $$B, h(B), h^2(B), h^3(B), h^4(B), h^5(B), h^6(B)=B$$
and $$C, h(C), h^2(C), h^3(C), h^4(C), h^5(C), h^6(C)=C$$
\subsubsection*{\bf Images of $d$}
Here, by contrast,  we will set $S_{1,d^6} =D$ (i.e. $S_{1,d^6} =:D$).

\begin{itemize}

\item
$S_{1,d}=d{\overline{d}}^{-1}=1$
\item

$S_{d,d}=dd{\overline{dd}}^{-1}= 1$
\item
$S_{d^d,d}=d^2d{\overline{d^3}}^{-1} =1$
\item
$S_{d^3,a}= d^3d{\overline{d^4}}^{-1} =1$
\item
$S_{d^4,d}= d^4d{\overline{d^4d}}^{-1}=1$
\item

$S_{d^5,d}= d^5d{\overline{d^5d}}^{-1}=d^5d{\overline{d^6d}}^{-1}= d^6\cdot 1= d^6 :=D$

\end{itemize}
We compute $\tau(d^jd^6d^{-j})= \tau(d^jDd^{-j})$ for all $j$ to see that, $h^j(D) = D.$

\subsubsection*{{\bf $\t(R)= \t([a,b][c,d])$} }
\smallskip
\noindent $\t(aba^{-1}b^{-1}c dc^{-1}d^{-1})=$
$$ S_{1,a} \cdot S_{{\overline{a}},b}\cdot S_{{\overline{aba^{-1}}}, a}^{-1} \\
 \cdot S_{{\overline{aba^{-1}b^{-1}}}, b}^{-1} \cdot S_{{\overline{aba^{-1}b^{-1}}}, c} \cdot
 S_{{\overline{aba^{-1}b^{-1}c} }, d}\cdot S_{ {\overline{aba^{-1}b^{-1}cdc^{-1}}}, c}^{-1}\cdot S^{-1}_{ {\overline{aba^{-1}b^{-1}cdc^{-1}d^{-1}}},  d}$$
Set $X = [A,B]$ and note that the $S_{*,d}$ words are all $1$. We obtain
\begin{equation}
XCh(C)=1
\end{equation}
\subsection*{{\bf Compute $\t(d^jRd^{-j})= \t(d^j[a,b][c,d]d^{-j})$}}
\vskip .05in
For $j = 1,2,3,4$ we obtain
\begin{equation}
h(X)h(C)h^2(C)=1
\end{equation}
\begin{equation}
h^2(X)h^2(C)h^3(C)=1
\end{equation}
\begin{equation}
h^3(X)h^3(C)h^4(C)=1
\end{equation}
\begin{equation}
h^4(X)h^4(C)h^5(C)=1
\end{equation}
But when we compute
\begin{equation} \t(d^5Rd^{-5})= \t(d^5[a,b][c,d]d^{-5}),\end{equation}
we obtain $D$ and $D^{-1}$ in the relation:
\begin{equation} \label{equation:**}
h^5(X)h^5(C)DC^{-1}D^{-1}=1
\end{equation}
That is, let $$Q =
S_{1,d} \cdot S_{d,d} \cdot S_{d^2,d} \cdot S_{d^3,d} \cdot S_{d^4}$$
We obtain
$$ Q \cdot  S_{d^5,a} \cdot S_{{\overline{d^5a}},b}\cdot _{{\overline{d^5aba^{-1}}}, a}^{-1}
\cdot S_{{\overline{d^5aba^{-1}b^{-1}}}, b}^{-1} \cdot
 S_{{\overline{d^5aba^{-1}b^{-1}}}, c} \cdot S_{{\overline{d^5aba^{-1}b^{-1}c} },d}
   \cdot S_{ {\overline{d^5aba^{-1}b^{-1}cdc^{-1}}}, c}^{-1}\cdot S^{-1}_{ {\overline{d^5aba^{-1}b^{-1}cdc^{-1}d^{-1}}},  d} \cdot Q^{-1}.$$
Since $Q=1$, we obtain equation \ref{equation:**}.
\smallskip

\subsubsection{\bf Eliminations}
\smallskip
Now we could eliminate images of $C$ from the six relations,
but instead we note that $XCh(C)^{-1}=1$ gives $h(C) = XC$ since $X$ is a commutator we have $h(C) \approx_h C$, where $\approx_h$ denotes homologous.

We end up with a single long relation  the product of $X$ and its five images under $h$ and the $CDC^{-1}D^{-1}$.

Since $g_0 =2$, the order of $G$ is six,  and there is no branching $2g-2=g(2g_0-2)$,  we have $2g =14$. Thus  we have the correct number of generators and their inverses  in the long  defining relation.

Now $X = [A,B]$ so that $h^j(A)$ and $h^j(B)$ are in the basis for all $j=1,...,6$.
We note that  $h^j(D) \approx_h D$ and $h^j(C) \approx_h C$ for all $j$. That is, we have twelve generators of type (1) and two generators of type (4). Types (1) and (4) refer to the types in definition \ref{defn:prep}.

Wirth respect to this basis, the matrix of the action of any element of $G$ will have two $6 \times 6$ blocks that are permutation matrices and one $2 \times 2$ identity matrix.

\section{The $G= \mathbb{Z}_2\times \mathbb{Z}_2$ example, example \#2} \label{section:extwo}

{\sl The Fuchsian group $\G_0$;}
$\G_0 = \langle a,b,c,d \;|\; abcd; a^2=b^2= c^2=d^2 \rangle $
\vskip .05in
{\sl The finite group:} $G= \mathbb{Z}_2 \times \mathbb{Z}_2= \langle g, h | g^2=h^2 =1, gh=hg \rangle$
\vskip .05in
{\sl The homomorphism $\phi$:} $\phi(a)=\phi(b) = g$, $\phi(c) = \phi(d) = h$, %$h^3 =1$
\vskip .05in {\sl The generating vector:} $(\phi(a), \phi(b), \phi(c), \phi(d) = (g,g,h,h)$
\vskip .05in {\sl
Coset representatives:} $1,b,d,bd$
\vskip .05in
{\sl Riemann Hurwitz:}
$2-2g = 4(2-2g_0) -4[(1-{\frac{1}{2}})+(1-{\frac{1}{2}}+(1-{\frac{1}{2}}+(1-{\frac{1}{2}}]$  since $g_0=0$
Thus $g=1$. The homology basis will have two curves.

\subsection{a-generators}

\smallskip

\vskip .05in

$\;$

$S_{1,a}= ab^{-1}=:A$

$S_{b,a}= ba$

$g(A) = \tau(bab^{-1}a^{-1})= S_{1,b}S_{b,a}S_{b,b}^{-1}S_{1,b}^{-1}= S_{b,a} = g(S_{1,a})$

$\tau(aa) = S_{1,a}S_{b,a}$

$\implies S_{b,a} = S_{1,a}^{-1}$

$S_{d,a} = da{\overline{bd}}^{-1}= dad^{-1}b^{-1}$

$(\tau(dab^{-1}d^{-1}) = S_{1,d}S_{d,a}S_{d,b}^{-1}S_{1,d}^{-1}= S_{d,a}S_{d,b}^{-1}$
-because all $S_{*,d} =1$, as we will see below.

$h(S_{1,a}) = S_{d,a}S_{d,b}^{-1}$

$h(S_{d,a})= \t(ddad^{-1}b^{-1}d^{-1}) =
S_{1,a}S_{b,d}^{-1}S_{d,b}^{-1}S_{1,d}^{-1} = S_{1,a}S_{b,d}^{-1}$.

$S_{bd,a}= bdad^{-1}$

$h(S_{bd,a}) = \t(dbdad^{-1})= S_{1,d}S_{d,b}S_{bd,d}S_{b,a}S_{1,d}= S_{d,b}S_{b,a}$

$S_{bd,a}= S_{d,a}^{-1}$

\subsection{b-generators}

\vskip .05in

$\;$

$S_{1,b} = bb^{-1} =1$

$S_{b,b} = b^2 \times 1 = 1$

$S_{d,b}= db(bd)^{-1}=dbd^{-1}b^{-1} = [d,b]\approx_h 1$

$S_{bd,b} = bdbd^{-1} = bdb^{-1}d^{-1} = [b,d] \approx_h  = 1$

$S_{bd,b}=S_{d,b}^{-1}$

Note $h(S_{d,b}) = \t(ddbd^{-1}b^{-1}d^{-1}) \t(bd^{-1}b^{-1}d^{-1} =  S_{1,b}S_{bd,d}^{-1}S_{d,b}^{-1}S_{1,d}^{-1} = S_{d,b}^{-1}=S_{bd,b}$.

$g(S_{d,b}) = \t(bdbd^{-1}b^{-1}b^{-1})= \t(bdbd^{-1}) =  S_{1,b}S_{b,d}S_{bd,b}S_{1,d}^{-1} = S_{bd,b}=S_{d,b}^{-1}$

$gh(S_{d,b}) = g(S_{b,d}^{-1}) =S_{d,b}$.

$g(S_{bd,b}) = \t(bbdbd^{-1}b^{-1}) = \t(dbd^{-1}b^{-1})= S_{b,d}=S_{bd,b}^{-1}$

$h(S_{bd,b}) = \t(dbdbd^{-1}d^{-1}) = \t(dbdb)= S_{d,b}=S_{bd,b}^{-1}$

We have computed all of the images of the $b$-words under the group.
\subsection{c-generators}

\vskip .05in
$\;$

$S_{1,c}=  cd^{-1}$

$S_{b,c} = bc(bd)^{-1}= bcd^{-1}b^{-1}$

$S_{d,c} = dc {\overline{(cd)}}^{-1} = dcd^{-2}  = dc$

$S_{bd,c}= bdc{\overline{bdc}}^{-1}= bdc(bd^2)^{-1} = bdcb^{-1}$

{\sl Images of the $c$ words under the group:}

$\;$

$g(S_{1,c}) = \t(bS_{1,c}b^{-1}) = \t(bcd^{-1}b^{-1} = S_{1,b}S_{b,c}S_{b,d}^{-1}S_{1,b}^{-1} = S_{b,c}$

$h(S_{1,c}) = \t(dS_{1,c}d^{-1}) = \t(dcd^{-1}d^{-1})=S_{d,c}$

$gh(S_{1,c}) = \t(bdcd^{-1}d^{-1}b^{-1}) = \t(bdcb^{-1}) = S_{1,b}S_{b,d}S_{bd,c}S_{1,b}^{-1} = S_{bd,c}$

We note $S_{b,c}=S_{bd,c}^{-1}$
and $S_{1,c} = S_{d,c}^{-1}$.
$g(S_{d,c})=  \t(bdcb^{-1}) = S_{1,b}S_{b,d}S_{bd,c}S_{1,b}^{-1}=S_{bd,c}$

$h(S_{d,c}) = \t(ddcd^{-1}) = S_{1,c}$

$gh(S_{d,c}) = g(S_{1,c}) = S_{b,c}$

$g(S_{b,c}) = \t(bbcd^{-1}b^{-1}b^{-1}) = \t(cd^{-1}) = S_{1,c}S_{1,d}^{-1} = S_{1,c}$

$h(S_{b,c}) = \t(dbcd^{-1}b^{-1}d^{-1}) = \t(dbcd^{-1}b^{-1}d^{-1}) =
S_{1,d}S_{d,b}S_{bd,c}S_{bd,d}^{-1}S_{d,b}^{-1}S_{1,d}^{-1} = S_{d,b}S_{db,c}S_{d,b}^{-1}
 $

$gh(S_{b,c}) =g(S_{d,b})g(S_{db,c})g(S_{d,b})^{-1})=S_{d,b}^{-1}g(S_{bd,c})S_{d,b}$

%$g(S_{db,c}) = \t(bdbcb^{-1}b^{-1}) = \t(bdbcb) = S_{1,b} S_{b,d} S_{bd,b}S_{d,c}S_{b,b}$

$g(S_{bd,c}) = \t(bbdc(bd^2){-1}b^{-1})= \t(bbdcd^{-1}d^{-1}b^{-1}b^{-1}) = S_{d,c}$

%**$g(S_{db,c}) = \t(bdbcb^{-1}b^{-1}) = \t(bdbc) = S_{1,b}S_{b,d}S_{bd,b}S_{d,c}= S_{bd,b}S_{d,c}$

\subsection{d-generators}

\vskip .05in

$\;$

$S_{1,d} = d \cdot d^{-1} =1$

$S_{b,d}= bd(bd)^{-1} =1$

$S_{d,d}= dd^{-1}=1$

$S_{bd,d}=bd^2(b)^{-1} =1$

\subsection{ $\t$ images}

\vskip .05in

$\;$

$\t(a^2) = S_{1,a}S_{b,a}=1$

$\t(ba^2b^{-1})= S_{1,b}S_{b,a}S_{1,a}S_{1,b}^{-1}= S_{b,a}S_{1,a}=1$

$\t(da^2d^{-1}) S_{1,d}S_{d,a}S_{bd,a}S_{1,d}^{-1} =S_{d,a}S_{bd,a}=1 $.

$\t(bda^2d^{-1}b^{-1})= S_{1,b} S_{b,d}S_{bd,a}S_{d,a}S_{b,d}{-1}S_{1,b}^{-1} =S_{bd,a}S_{d,a}=1 $.

\subsection{$\tau(Kb^2K^{-1})$}

\vskip .05in

$\;$
$\t(b^2) = S_{1,b}S_{b,b}=1$

$\t(bb^2b^{-1})= 1$

$\t(db^2d^{-1})=  S_{1,d}S_{d,b}S_{bd,b}S_{1,d}^{-1} =S_{d,b}S_{bd,b}=1 $.

$\t(bdb^2d^{-1}b^{-1})= S_{1,b} S_{b,d}S_{bd,b}S_{d,b}S_{*,d}{-1}S_{1,b}^{-1} =S_{bd,b}S_{d,b}=1 $.

\subsection{$\t(Kc^2K^{-1})$}

\vskip .05in

$\;$
$\t(c^2) = S_{1,c}S_{b,c}=1$

$\t(bc^2b^{-1})= S_{1,b}S_{b,c}S_{bd,c}S_{1,b}^{-1}=S_{b,c}S_{bd,c}$

$\t(dc^2d^{-1})=  S_{1,d}S_{d,c}S_{1,c}S_{1,d}^{-1} =S_{d,c}S_{1,c}$.

$\t(bdc^2d^{-1}b^{-1})= S_{1,b} S_{b,d}S_{bd,c}S_{d,c}S_{*,d}^{-1}S_{1,b}^{-1} =S_{bd,c}S_{d,c}=1 $.

\subsection{$\t(Kd^2K^{-1})=1 \forall K$}

\vskip .05in
$\;$

\subsection{$\t(KRK^{-1})$} where $R= abcd$
%\vskip .05in
\begin{enumerate}

\item $\t(abcd) = S_{1,a}S_{{\overline{a}},b}S_{{\overline{ab}},c}S_{{\overline{abc}},d}
    = S_{1,a}S_{b,b}S_{1, c}S_{d,d} = S_{1,a}S_{1,c}=1 $

\bigskip
\item  $\t(babcdb^{-1})= S_{1,b}S_{b,a}S_{{\overline{ba}},b}S_{{\overline{bab}},c}S_{{\overline{babc}},d}S_{1,b}^{-1} = S_{b,a}S_{b,c}=1$
\bigskip

\item  \label{item:QRT}$\t(dabcdd^{-1})= S_{1,d} S_{d,a}S_{{\overline{da}},b}S_{{\overline{dab}},c}S_{{\overline{dabc}},d}S_{1,d}^{-1}=
S_{d,a}S_{bd,b}S_{b,c}=1$
\bigskip

\item  \label{item:QinvRinvTinv}$\t(bdabcdd^{-1}b^{-1}) =
S_{1,b}S_{b,d}S_{bd,a}S_{{\overline{bda}},b}S_{{\overline{bdab}},c}
S_{{\overline{bdabc}},b}^{-1} = S_{bd,a}S_{d,b}S_{bd,c}=1$

%\section{Combining}
\end{enumerate}

\subsection{Combining}

\vskip .05in

$\;$

Let $M=:S_{d,a}$ so that $M^{-1} =: S_{bd,a}$ and let $T =: S_{b,c}$ so that $T^{-1} =: S_{bd,c}$.
Let $Z = S_{b,d}$. Since $S_{d,b} = S_{bd,b}^{-1}$, we have by \ref{item:QRT}  $MZT=1$ and by \ref{item:QinvRinvTinv} $M^{-1}Z^{-1}T^{-1}=1$ so that $Z^{-1} = MT$ or $Z = T^{-1}M^{-1}$, yielding $$MTM^{-1}T^{-1} = 1.$$
We note:
$$g(M) = \t(bS_{d,a}b^{-1}) = \t(bdad^{-1}b^{-1}b^{-1}) = S_{1,b}S_{b,d}S_{bd,a}S_{1,d}^{-1} = S_{bd,a} = S_{d,a}^{-1}=M^{-1}.$$
$$h(M) = \t(d S_{d,a}d^{-1})
= \t(ddad^{-1}b^{-1}d^{-1})= S_{1,a}= S_{b,a}^{-1}=M^{-1}$$
$$gh(M) =M$$
$$g(T) = \t(bS_{b,c}b^{-1}) = \t(bbcd^{-1}b^{-1}b^{-1}) = \t(cd^{-1}) =$$
$$g(T) = S_{1,c}S_{1,d}^{-1}= S_{1,c} = S_{b,c}^{-1} =T^{-1}.$$
$$h(T) = \t(dS_{b,c}d^{-1}) = \t(dbcd^{-1}b^{-1}d^{-1}) = S_{1,d}S_{b,d}S_{bd,c}S_{b,d}^{-1}S_{b,b}^{-1}S_{1,d}^{-1}= S_{bd,c}=S_{d,c}^{-1} = T^{-1}$$
$$ \mbox{ and } gh(T) = T.$$

We are left with two generators $M$
 and $T$ each of type (4) definition \ref{defn:prep}.  Finally we can double check that the elements in our group fix the correct number of points on the surface using the fact that the number of fixed points of an element is the same as $2$ - the trace of its action on homology by the Lefshetz fixed point formula \cite{MacB}.

\section{Example \#3: different $\G_0$, but the same $G = \ZZ_2
\times \ZZ_3$} \label{section:exthree}
\vskip .05in

$\;$

{\sl The Fuchsian group $\G_0$:}
$\G_0 = \langle a,b,c,d,e \;|\; abcde; a^2=b^2; c^3=d^3=e^3 \rangle $
\vskip .05in
{\sl The finite group $G$:}
$G= \mathbb{Z}_2 \times \mathbb{Z}_3 = \langle g, h | g^2=h^3 =1, gh=hg \rangle$
\vskip .05in
{\sl The surface kernel map $\phi$:} $\phi(a)=\phi(b) = g$, $\phi(c) = \phi(d) =\phi(e) =h$, %$h^3 =1$
\vskip .05in
{\sl The generating vector:} $(\phi(a), \phi(b), \phi(c), \phi(d), \phi(e)) = (g,g,h,h,h)$
\vskip .05in
{\sl Coset representatives:} $1,b,d,d^2,bd, bd^2$

\vskip .05in

Use the Riemann-Hurwitz theorem to compute that since $g_0=0$,  $g=4$, so that a homology basis will have $8$ elements.

\subsection{a-generators}

\vskip .05in
$\;$

$S_{1,a}= ab^{-1}=:A$

$S_{b,a}= ba$

$\tau(bab^{-1}a^{-1})= S_{1,b}S_{b,a}S_{b,b}^{-1}S_{1,b}^{-1}= S_{b,a}.$
  Later we will write $S_{b,a} = g(S_{1,a}).$

$\tau(aa) = S_{1,a}S_{b,a}$

  $\implies S_{b,a} = S_{1,a}^{-1}$

$S_{d,a} = da{\overline{bd}}^{-1}$

$\tau(dab^{-1}d^{-1}0 = S_{1,d}S_{d,a}S_{d,d}^{-1}S_{1,d}^{-1}= S_{d,a}.$ Later we will write  $h(S_{1,a}) = S_{d,a}.$

$S_{d^2,a}= d^2a(bd^2)^{-1} = d^2ad^{-2}b^{-1}= [d^2,a]ab^{-1}.$

 Note $S_{d^2,a} = \approx_h S_{1,a}$

$S_{bd,a}$

$S_{bd^2,a}$

\subsection{b-generators}

$\;$

\vskip .05in

$S_{1,b} = bb^{-1} =1$

$S_{b,b} = b^2$

Since this is a relator. We can  remove $S_{b,b}$ from the generating set if we replace each of its occurrences by $1$.

$S_{d,b}= db(bd)^{-1}=dbd^{-1}b^{-1} = [b,d].$   Thus  $S_{d,b} \approx_h 1$.

$S_{d^2,b} = d^2b(bd^2)^{-1} = [b,d^2].$ Thus $S_{d^2,b} \approx_h 1.$

$S_{bd,b} = bdbd^{-1} = b[d,b]b.$ Thus $S_{bd,b} \approx_h 1.$

$S_{bd^2,b}= bd^2bd^-2= [b,d^2] = b[d^2,b]b.$  Thus  $S_{bd^2,b} \approx_h 1.$

\subsection{c-generators}

$\;$

\vskip .05in
$S_{1,c}=
 cd^{-1}$

$S_{b,c} = bc(bd)^{-1}= bcd^{-1}b^{-1}.$ Since $S_{b,c} = [b,c]cb^{-1}$, $S_{b,c}  \approx_h 1  \approx cb^{-1}.$
and
$\t(S_{b,c})=  \t(bbc{\overline{bc}^{-1}b^{-1}} = \t(bbcd^{-1}b^{-1}b^{-1}) \t(cd^{-1}= S{1,c}S{1,d}^{-1} = S_{1,c}.$

Thus
$S_{d,c} = dc (cd)^{-1} = dcd^{-2}$

Thus $S_{d,c}= dcd^{-1}c^{-1}cd^{-1}= [d,c^{-1}]cd^{-1} =\approx_h 1 cd^{-1} \approx_h cd^{-1}= S_{1,c}^{-1}$

$S_{d^2,c} = d^2c {\overline{d^2c}}^{-1} = d^2c$

and

$h(S_{d^2,c}) \t(dd^2cd^{-1})=\t(cd^{-2} =S_{1,c}S_{1,d}^{-1}= S_{1,c}.$

$S_{bd,c}= bdc{\overline{bdc}^{-1}}= bdc(bd^2)^{-1} = bdcd^{-2}b^{-1}= bdcd^{-1}b^{-1} =$

$b[d,c]cd^{-1}b^{-1}= [b,[d,c]][d,c]bcd^{-1}b^{-1}$

$=[b,[d,c]][d,c]bcb^{-1}c^{-1}cbd^{-1}b^{-1}$

$=[b,[d,c]][d,c] cd^{-1}dbd^{-1}b^{-1} b^{-1}= [b,[d,c]][d,c]cd^{-1}[d,b].$

Thus  $S_{bd,c} \approx_h   cd^{-1}$

$S_{bd^2,c} = bd^2cb^{-1}= bd^{-1}cb^{-1}$

$S_{bd^2,c}= bd^{-1}cd c^{-1}cd^{-1}b^{-1} = b^{-1}[d^{-1},c] cd^{-1}b$

$ S_{bd^2,c} \approx_h b^{-1}cd^{-1}b = (cd^{-1})[cd^{-1},b^{-1}]$

$h(S_{bd^2,c}) = [b,d]S_{b,c}[b,d]^{-1}$

 $h(S_{b,c}) \approx_h  \t(dcb^{-1}d^{-1}= \approx_h S_{1,d}S_{d_2,c}S_{bd^2,b}^{-1}S_{bd,d}^{-1} \approx_h S_{d^2,c} $

\subsection{d-generators}

\vskip .05in

$\;$

$S_{1,d} = d d^{-1} =1$

$S_{b,d}= bd(bd)^{-1} =1$

$S_{d,d}= dd^{-1}=1$

$S_{d^2,d}= d^3 = 1$

$S_{bd,d}=bd^2(bd^2)^{-1} =1$

$S_{bd^2,d} = bd^dd b^{-1}= bd^3b^{-1}$

\subsection{e-generators}

\vskip .05in

$\;$

$S_{1,e}= ed^{-1}$

$S_{b,e} = be(bd)^{-1}= bed^{-1}b^{-1}.$

\vskip .02in
Thus  $g(S_{1,e}) = \t(bS_{1,b}b^{-1}) =  \t(bed^{-1}b^{-1}) =S_{b,e}$

$S_{d,e} = de (ed)^{-1}= ded^{-2} $

Thus $h(S_{1,e}) = \t(ded^{-2})=
 S_{1,d}S_{d,e}S_{d,d}^{-1}S_{1,d}^{-1}= S_{d,e}, $

 and we have also: $S_{d,e} =ded^{-1}e^{-1}ed^{-1}) =([d,e]ed^{-1}) \approx_h S_{1,e}$

 $S_{d^2,e} = d^2e= d^{-1}e$

Thus
 $h(S_{d,e}) = \t(dS_{d,e}d^{-1})  = \t(dded^{-2}d^{-1} = \t(dde) =S_{1,d}S_{d,d}S_{d^2,e}=S_{d^2,e}.$

$h(S_{d^2,e}) = \t(ddde{\overline{ d^2e}}^{-1}d^{-1})=
\t(ed^{-1})= S_{1,d}S_{1,d}^{-1} =S_{1,e}.$

$S_{bd,e}= bde{\overline{bde}^{-1}}= bde(bd^2){-1} = bded^{-2}b^{-1}.$

Thus $S_{bd,e}  = bded^{-1}e^{-1}ed^{-1}b^{-1} = bd[e,d^{-1}]ed^{-1}b^{-1}
\approx_h bded^{-1}b^{-1} = bded^{-2}db^{-1}= bd \cdot ed^{-1}e^{-1}d \cdot d^{-1}e  \cdot d^{-1}b^{-1} = bd [e,d^{-1}] d^{-1}e d^{-1}b^{-1} = \mbox{product of some commutators }d^{-1}e \approx_h d^{-1}e = S_{d^2,e}. $

$S_{bd^2,e} = bd^2eb^{-1} = bd^{-1}eb^{-1}$

$S_{bd^2,e} = bd^{-1} ede^{-1}ed^{-1}  b^{-1}= b[d^{-1},e]ed^{-1}b^{-1} = b[d^{-1},e]S_{1,e}b^{-1}= [b,[d^{-1}e,e]S_{1,e}]S_{1,e} \approx_h S_{1,e} $
%g_b(S_bd^2, e} \approx_h d^2e$

\subsection{ $\t$ images}

\vskip .05in

$\;$

$\t(a^2) = S_{1,a}S_{b,a}=1$

$\t(ba^2b^{-1})= S_{1,b}S_{b,a}S_{1,a}S_{1,b}^{-1}= S_{b,a}S_{1,a}=1$

$\t(da^2d^{-1}) S_{1,d}S_{d,a}S_{bd,a}S_{1,d}{-1} =S_{d,a}S_{bd,a}=1 $.

$\t(d^2a^2d^{-2})= S_{1,d}S_{d,d}S_{d^2,a}S_{bd^2,a}S_{d,d}^{-1}S_{1,d}^{-1} =S_{d^2,a}S_{bd^2,a}=1$

$\t(bda^2d^{-1}b^{-1})= S_{1,b} S_{b,d}S_{bd,a}S_{d,a}S_{b,d}{-1}S_{1,b}^{-1} =S_{bd,a}S_{d,a}=1 $.

$\t(bd^2a^2d^{-2}b^{-1})= S_{1,b}^{-1}S_{b,d}S_{bd,d}S_{bd^2,a}S_{d^2,a}S_{bd,d}^{-1}S_{b,d}{-1}S_{1,b}^{-1} =
S_{bd^2,a}S_{d^2,a}=1 $

Note we can thus eliminate $a-$ words at some point and combine the 6 relators into 3.

\subsection{$\t(Kb^2K^{-1})$}

\vskip .05in

$\;$

$\t(bb) = S_{1,b}S_{b,b} = 1$
but each these two words are already the identity, thus the relator is $S_{b,b}=1$

$\t(bbbb^{-1})$ yields nothing that is, $\t(b^2) = t(bbbb^{-1})=1 $.

$\t(dbbd^{-1}) = S_{1,d}S_{d,b}S_{bd,b}S_{1,d}^{-1} \implies S_{d,b}=S_{bd,b}^{-1}$

$\t(bdbbd^{-1}b^{-1})= S_{1,b}S_{b,d}S_{bd,b}S_{d,b}S_{b,d}^{-1}S_{1,b} \implies
S_{d,b}= S_{bd,b}^{-1}$

$\t(ddbbd^{-1})= S_{1,d}S_{d,d}S_{d^2,b}S_{bd^2,b}S_{d,d}^{-1} =1 \implies
S_{d^2,b} = S_{bd^2,b}^{-1}$.

\subsection{$\t(KcccK^{-1})$}
\vskip .05in

$\;$

$\t(ccc) = S_{1,c}S_{d,c}S_{d^2,c} = 1$

I.e. $Ch(C)h^2(C) =1$.

$\t(bcccb^{-1})=1 \implies g(C)gh(C) gh^2(C) =1$.

\subsection{$\t(KeeeK^{-1})$}
\vskip .05in
$\t(eee) = S_{1,e}S_{d,e}S_{d^2,e} = 1$

Eventually we will have,  $Eh(E)h^2(E) =1$ if we set $E=S_{1,e}$ and

$\t(beeeb^{-1})=1 \implies g(E)gh(E) gh^2(E) =1$.

\subsection{$\t(ddd)$ }

\vskip .05in
$\;$

$$\t(ddd) = S_{1,d}S_{d,d}S_{d^2,d}$$

$$\t(bdddb^{-1})= S_{1,b}S_{b,d}S_{bd,d}S_{bd^2,d}S_{1,b}^{-1} =
S_{b,d}S_{bd,d}S_{bd^2,d}$$
%\end{document}

\subsection{$\t(KRK^{-1})$}  where $R= abcde.$

\begin{enumerate}

\item $\t(abcde) =S_{1,a}S_{{\overline{a}},b}S_{{\overline{ab}},c}S_{{\overline{abc}},d}
    S_{{\overline{abcd}},e}.$

{\sl Note that all of the $S_{*,d}$ are $1$. Let $A=S_{1,a}, C= S_{1,c}, \mbox{ and }  E= S_{1,e}$.}

$$\t(abcde) = A \cdot S_{b,b}\cdot C \cdot S_{d,d} \cdot h^2(E) =1.$$

{\sl Continue computing.}
%\begin{enumerate}
\item $\t(babcdeb^{-1})=1$ gives
$$S_{1,b} S_{b,a} S_{{\overline{ba}},b}S_{{\overline{bab}},c}S_{{\overline{babc}},d}
    S_{{\overline{babcd}},e}S_{1,b}^{-1}$$
    \item $\t(dabcded^{-1}) = 1$ gives
$$S_{1,d} S_{d,a} S_{{\overline{da}},b}S_{{\overline{dab}},c}S_{{\overline{dabc}},d}
    S_{{\overline{dabcd}},e}S_{1,d}^{-1}$$
\item $\t(d^2abcded^{-2}) = 1$
$$S_{1,d}S_{d,d}S_{d^2,a} S_{{\overline{d^2a}},b}S_{{\overline{d^2ab}},c}S_{{\overline{d^2abc}},d}
    S_{{\overline{d^2abcd}},e}S_{d,d}^{-1}S_{1,d}^{-1}$$
\item $\t(bdabcded^{-1}b{-1}) = 1$ gives
$$S_{1,b}S_{b,d} S_{bd,a} S_{{\overline{bda}},b}S_{{\overline{bdab}},c}S_{{\overline{bdabc}},d}
    S_{{\overline{bdabcd}},e}S_{b,d}^{-1}S_{1,b}^{-1}$$
\item $\t(bd^2abcded^{-2}b^{-1} ) = 1$ gives
$$S_{1,b}S_{b,d}S_{bd,d} S_{bd^2,a} S_{{\overline{bd^2a}},b}S_{{\overline{bd^2ab}},c}S_{{\overline{bd^2abc}},d}
    S_{{\overline{bd^2abcd}},e}S_{bd,d}^{-1}S_{b,d}^{-1}S_{1,b}^{-1}$$
\end{enumerate}

We solve the above relations and observe
\begin{enumerate}

\item $\t(abcde) = 1$ gives
    $$S_{1,a}^{-1} = S_{{\overline{a}},b}S_{{\overline{ab}},c}S_{{\overline{abc}},d}
    S_{{\overline{abcd}},e}$$
\item $\t(babcdeb^{-1})=1$ gives (since $S_{1,b} =1$),
$$S_{b,a}^{-1} =  S_{{\overline{ba}},b}S_{{\overline{bab}},c}S_{{\overline{babc}},d}
    S_{{\overline{babcd}},e}$$
\item $\t(dabcded^{-1}) = 1$ gives (since $S_{1,d}=1$)
$$ S_{d,a}^{-1} = S_{{\overline{da}},b}S_{{\overline{dab}},c}S_{{\overline{dabc}},d}
    S_{{\overline{dabcd}},e}$$
\item $\t(d^2abcded^{-2}) = 1$ gives (since $S_{1,d}=1$ and $S_{d,d}=1$,)
$$S_{d^2,a}^{-1} =  S_{{\overline{d^2a}},b}S_{{\overline{d^2ab}},c}S_{{\overline{d^2abc}},d}
    S_{{\overline{d^2abcd}},e}$$
\item $\t(bdabcded^{-1}b{-1}) = 1$ gives
(since $S_{1,b}=1$ and $S_{b,d}=1 $),
$$S_{bd,a}^{-1} =  S_{{\overline{bda}},b}S_{{\overline{bdab}},c}S_{{\overline{bdabc}},d}
    S_{{\overline{bdabcd}},e}$$
\item $\t(bd^2abcded^{-2}b^{-1} ) = 1$ gives
(since $S_{1,b}=1$, $S_{b,d}=1$ and $S_{bd,d}=1$),
$$S_{bd^2,a}^{-1} =  S_{{\overline{bd^2a}},b}S_{{\overline{bd^2ab}},c}S_{{\overline{bd^2abc}},d}
    S_{{\overline{bd^2abcd}},e}$$
\end{enumerate}
\subsection{\bf Elimination to three relations:}

\vskip .05in

$\;$

Now we have relators $$S_{1,a}S_{b,a}=1$$
$$S_{bd,a}S_{d,a}=1$$
$$S_{bd^2,a}S_{d^2,a}=1$$

We can, therefore, eliminate the $\t(Ka^2K^{-1})=1$ relations for all $K$.
and obtain three  relators from the six long relators and eliminate the six $a^2$-relators and the six $a-$generators.

\begin{enumerate}
\item $$ S_{{\overline{ba}},b}S_{{\overline{bab}},c}S_{{\overline{babc}},d}
    S_{{\overline{babcd}},e} S_{{\overline{a}},b}S_{{\overline{ab}},c}S_{{\overline{abc}},d}
    S_{{\overline{abcd}},e}$$
\item $$S_{{\overline{bda}},b}S_{{\overline{bdab}},c}S_{{\overline{bdabc}},d}
    S_{{\overline{bdabcd}},e}    S_{{\overline{da}},b}S_{{\overline{dab}},c}S_{{\overline{dabc}},d}
    S_{{\overline{dabcd}},e}$$
\item $$ S_{{\overline{bd^2a}},b}S_{{\overline{bd^2ab}},c}S_{{\overline{bd^2abc}},d}
    S_{{\overline{bd^2abcd}},e}S_{{\overline{d^2a}},b}S_{{\overline{d^2ab}},c}S_{{\overline{d^2abc}},d}
    S_{{\overline{d^2abcd}},e}$$
\end{enumerate}

Use the fact that $S_{1,d}S_{d,d}S_{d^2,d}$ is a realtor to combine these three relators into one relator and eliminate the generators $S_{1,d}$, $S_{d,d}$ and $S_{d^,d}$, but
First rewrite without the ${\overline{x}}$
 \begin{enumerate}
 \item $S_{1,b}S_{b,c}S_{bd,d}S_{bd^2,e}S_{b,b}S_{1,c}S_{d,d}S_{d^2,e}$

 \item $S_{d,b}S_{bd,c}S_{bd^2,d}S_{b,e}S_{bd,b}S_{d,c}S_{d^2,d}S_{1,e}$

 \item $S_{d^2,b} S_{bd^2,c} S_{b  , d} S_{bd  ,e}S_{bd^2,b}S_{d^2,c} S_{1 ,d} S_{d,e}$
 \end{enumerate}
 and then
 \begin{enumerate}
 \item $S_{d,d}^{-1}  = S_{d^2,e}S_{1,b}S_{b,c}S_{bd,d}S_{bd^2,e}S_{b,b}S_{1,c}$
 \item $S_{d^2,d}^{-1} = S_{1,e}S_{d,b}S_{bd,c}S_{bd^2,d}S_{b,e}S_{bd,b}S_{d,c}$
 \item $S_{1,d}^{-1} = S_{d,e}S_{d^2,b} S_{bd^2,c} S_{b  , d} S_{bd  ,e}S_{bd^2,b}S_{d^2,c} $

 \end{enumerate}

We obtain  $$S_{d^d,d}^{-1}S_{d,d}^{-1}S_{1,d}^{-1} =$$
$$S_{d^2,e}S_{1,b}S_{b,c}S_{bd,d}S_{bd^2,e}S_{b,b}S_{1,c}
S_{1,e}S_{d,b}S_{bd,c}S_{bd^2,d}S_{b,e}S_{bd,b}S_{d,c}
S_{d,e}S_{d^2,b} S_{bd^2,c} S_{b  , d} S_{bd  ,e}S_{bd^2,b}S_{d^2,c}$$

Finally, replace $S_{d^2,e}$ by $(S_{1,e}S_{d,e})^{-1}$, $S_{bd^2,e}$ by $(S_{b,e}S_{bd,e})^{-1}$,
$S_{d^2,c}$ by $(S_{1,c}S_{d,c})^{-1}$, $S_{bd^2,c}$ by $(S_{b,c}S_{bd,c})^{-1}.$

Eliminate the $b-$ generators and relations, to obtain one single relation (displayed on two lines, but a product):

$$(S_{1,e}S_{d,e})^{-1}S_{b,c}S_{bd,d}(S_{b,e}S_{bd,e})^{-1}S_{1,c}
S_{1,e}S_{bd,c}S_{b,e}S_{d,c}
S_{d,e}(S_{b,c}S_{bd,c})^{-1}S_{bd  ,e}(S_{1,c}S_{d,c})^{-1} \times $$

$$S_{d^2,e}S_{b,c}S_{bd,d}S_{bd^2,e}S_{1,c}
S_{1,e}S_{bd,c}S_{b,e}S_{d,c}
S_{d,e}S_{bd^2,c} S_{bd  ,e}S_{d^2,c}$$

And lastly eliminate the remaining $d-$generators and relations to obtain
one long relation involving eight generators for the fundamental group of $\G$ and their inverses
$$(S_{1,e}S_{d,e})^{-1}S_{b,c}
(S_{b,e}S_{bd,e})^{-1}S_{1,c}
S_{1,e}
S_{bd,c}
S_{b,e}S_{d,c}
S_{d,e}
(S_{b,c}S_{bd,c})^{-1}
 S_{bd  ,e}
(S_{1,c}S_{d,c})^{-1}$$

\subsection{Final homology computation}

\vskip .05in

$\;$

Recall that  $g$ and $h$ denote the action of elements of $\G_0$ on $\G$ by conjugation.

Now we may assume that  $E$ is the image of $S_{1,e}$ in the first homology group
so that $h(E), h^2(E)$ denote the images under $h$ with $g(E), hg(E)=gh(E)$, $h^2g(E) = gh^2(E)$ the other images of $E$.

Similarly let $C$ be the image in the first homology group of $S_{1,c}$.

We have a set of generators for homology of the form:
$$E, h(E),g(E), gh^2(E),gh(E), C, h(C),g(C), gh^2(C),gh(C)$$

We turn this into a homology basis with

$E$ and $g(E)$ are in the basis and
$C$ and $g(C)$ are in the basis.
The same holds for $E^{-1}$ and $C^{-1}$ and their images under $g$.

$E$ and $h(E)$ are in the basis with $h^2(E) \approx_h -(E + h(E))$
$g(E)$ and $gh(E)$ are in the basis with $gh^2(E) \approx_h -(g(E) + h^2g(E))$

$C$ and $h(C)$ are in the basis with $h^2(C) \approx_h -(C + h(C))$

$g(C)$ and $gh(C)$ are in the basis with $gh^2(C) \approx_h -(g(C) + hg(C))$

with appropriate images relations for the curves under the other generators for $\G_0/\G$ and $g(C)$
\section{The fixed points and matrices}

We compute the matrices of the action on homology with respect the {\sl ordered}  basis
$E, g(E), h(E), gh(E), C, g(C), h(C), gh(C)$
\vskip .05in

$\;$

Let $P= \left(
  \begin{array}{cccc}
    0 & 1 & 0 & 0 \\
    1 & 0 & 0  & 0 \\
    0 & 0 & 0 & 1 \\
    0 & 0 & 1 & 0 \\
  \end{array}
\right), Q= \left(
  \begin{array}{cccc}
    0 & 0 & 1 & 0 \\
    0 & 0 & 0  & 1 \\
    -1 & 0 & -1 & 0 \\
    0 & -1 & 0 & -1
     \\
  \end{array}
\right), \mbox{and} \;  T= \left(
  \begin{array}{cccc}
    0 & 0 & 0 & 1 \\
    0 & 0 & 1  & 0 \\
    0 & -1 & 0 & -1 \\
   -1 & 0 & -1 & 0 \\
  \end{array}
\right)$

\vskip .05in

$\;$

The matrices of $g, h$ and $gh$ with respect to  the basis $$E,g(E),h(E), gh(E),C,g(C),h(C), gh(C)$$
break up into blocks where $0_n$ represents the $n \times n$
the matrix of all zeros for an integer $n$ and thus $0_4$ the $4 \times 2$ such matrix.

\vskip .05in

$\;$

$M_g = \left(
  \begin{array}{cccc}
    P & 0_4 \\
    0_4 & P \\
      \end{array}
\right), M_h = \left(
  \begin{array}{cccc}
    Q & 0_4 \\
    0_4 & Q \\

  \end{array}
\right), \mbox{and} \; M_{gh} = \left(
  \begin{array}{cccc}
    T & 0_4 \\
    0_4 & T \\

  \end{array}
\right)$
%\end{document}

\vskip .05in
Using the Lefschetz fixed point formula \cite{MacB} we can double check that these transformations have the correct number of fixed points (as defined by the original surface kernel map) by computing $2-\Tr M_h = 6$, $2-\Tr M_g =2$, and  $2-\Tr M_{gh} =0$ where $\Tr$ denotes the trace of a matrix.
\vskip .01in

Observe that $M_g$ has four non-zero blocks that are $2 \times 2$ permutation matrices.
The matrix $M_h$ is built upon an element of type (2) and if the ordered basis is rearranged to be
$E,h(E), g(E), gh(E), C, h(C), g(C), gh(C)$ we have if
\smallskip

 $Q_0= \left(
  \begin{array}{cccc}
    0 & 1 & 0 & 0 \\
    -1 & -1 & 0  & 0 \\
    01 & 0 & 0 & 1 \\
    0 & 01 & -1 & -1
     \\
  \end{array}
  \right)$,  then $h$ has the matrix   $\left(
  \begin{array}{cccc}
    Q_0 & 0_4 \\
    0_4 & Q_0 \\

  \end{array}
\right)$.

If the ordered basis is $E, h^2(E), gh(E), g(E), C, h^2(C), gh(C), g(C)$, let
\smallskip

$T_0= \left(
  \begin{array}{cccc}
    0 & 0 & 1 & 0 \\
    0 & 0 & 0  & 1 \\
    0 & 1 & 0 & 0 \\
   -1 & -1& 0 & 0 \\
  \end{array}
\right)$
so
that $gh$ then has the matrix
 $\left(
  \begin{array}{cccc}
    T_0 & 0_4 \\
    0_4 & T_0 \\

  \end{array}
\right)$.
\vskip .5in

{\bf Acknowledgement} The author thanks Marston Conder, Gareth Jones and David Singerman for some useful conversations about connections between Schreier-Reidemeister Theory and conformal automorphism groups. The author also thanks the referee for some useful suggestions.

\bibliographystyle{amsalpha}

\end{document}